\documentclass{article}
\usepackage{amsmath, amsthm, amsfonts, amssymb}
\usepackage{hyperref}
\usepackage[latin1]{inputenc}

\makeatletter \@addtoreset{equation}{section}

\makeatletter

\def\beg   {\begin{theorem}}   \def\ee   {\end{theorem}}
\def\be   {\begin{equation}}   \def\ee   {\end{equation}}
\def\ba   {\begin{array}}      \def\ea   {\end{array}}
\def\bea  {\begin{eqnarray}}   \def\eea  {\end{eqnarray}}
\def\bean {\begin{eqnarray*}}  \def\eean {\end{eqnarray*}}

\newtheorem{lemma}{Lemma}[section]
\newtheorem{theorem} [lemma]{Theorem}

\newtheorem{definition}[lemma] {Definition}

\newcommand{\mC}{\ensuremath{\mathbb{C}}}

\begin{document}

\vspace{4cm}
\begin{center} \LARGE{\textbf{ Uniqueness of differential polynomials sharing one value  }}
 \end{center}
 \vspace{1cm}
 \begin{center} \bf{Banarsi Lal$^{1 }$, \quad Kuldeep Singh Charak$^{2}$ }
\end{center}

\begin{center}
 Department of Mathematics, University of Jammu,
Jammu-180 006, INDIA.\\
{$^{1}$ E-mail: banarsiverma644@gmail.com }\\

{$^{2}$ E-mail: kscharak7@rediffmail.com }
\end{center}

\bigskip
\begin{abstract}
 We prove some uniqueness results which improve and generalize results of Jiang-Tao Li and Ping Li[{\it Uniqueness of entire functions concerning differential polynomials. Commun. Korean Math. Soc. 30 (2015), No. 2, pp. 93-101}].
\end{abstract}

\vspace{1cm}\noindent \textbf{Keywords: } Meromorphic functions, sharing of values, differential polynomials, Nevanlinna theory.

\vspace{0.5cm} \noindent\textbf{AMS subject classification: 30D35, 30D30}

\vspace{4cm}

\normalsize
\newpage

\section{Introduction}
Let $f$ be a non-constant meromorphic function in the complex plane $\mC$. We assume that the reader is familiar with the standard notions of the Nevanlinna value distribution theory such as $T(r, f), \; m(r, f),\;N(r, f)$  (see e.g., \cite{HAY}). \\

For $a \in  \mC \cup \{\infty \}$, we say that two meromorphic functions $f$ and $g$ share a CM, if $f - a$ and $g - a$ have the same set of zeros with same multiplicities, and if we do not consider the multiplicities then $f$ and $g$ are said to share a IM.

\medskip

In \cite{CM'}, C.C. Yang posed the following question:\\

{\bf Question:} What can be said about two entire functions $f$ and $g$, when they share $0$ CM and their derivatives share $1$ CM $?$

\medskip

In 1990, Yi \cite{HX, HX'}, answered the above question by proving: {\it Let $f$ and $g$ be two non-constant entire functions such that $f$ and $g$ share $0$ CM. If $f^{(k)}$ and $g^{(k)}$ share the value $1$ CM and $\delta (0, f) > 1/2$, where $k$ is non-negative integer, then $f \equiv g$ unless $f^{(k)}.g^{(k)} \equiv 1$;} and for meromorphic functions he proved:  {\it Let $f$ and $g$ be two non-constant meromorphic functions such that $f$ and $g$ share $0$ and $\infty$ CM. If $f^{(k)}$ and $g^{(k)}$ share the value $1$ CM and $2\delta (0, f) + (k + 2)\Theta(\infty , f) > k + 3$, where $k$ is non-negative integer, then $f \equiv g$ unless $f^{(k)}.g^{(k)} \equiv 1$.}  

\medskip

For a non-constant meromorphic function $h$, we denote by 
$$L(h) = h^{(k)} + a_1 h^{(k-1)} + a_2 h^{(k-2)} + ... + a_{k-1} h' + a_k h,$$
the differential polynomial of $h$, where $a_1,\;a_2,\;...,\;a_k$ are finite complex numbers and $k$ is a positive integer. We denote the order and lower order of $h$ by $\lambda(h)$ and $\mu(h)$, respectively. Also by $\sigma(h)$ and $\sigma(1/h)$, we denote the exponent of convergence of zeros and poles of $h$ respectively.\\

Recently, Jiang-Tao Li and Ping Li $\cite{LL}$ generalized first result of Yi(as stated above) for entire fuctions as

\medskip

{\bf Theorem A.}
Let $f$ and $g$ be two non-constant entire functions such that $f$ and $g$ share $0$ CM. Suppose $L(f)$  and $L(g)$ share $1$ CM and $\delta(0, f) > 1/2$.
If $\lambda (f) \neq 1$, then $f \equiv g$ unless $L(f).L(g) \equiv 1$.

\medskip

{\bf Theorem B.}
Let $f$ and $g$ be two non-constant entire functions such that $f$ and $g$ share $0$ CM. Suppose $L(f)$  and $L(g)$ share $1$ IM and $\delta(0, f) > 4/5$.
If $\lambda (f) \neq 1$, then $f \equiv g$ unless $L(f).L(g) \equiv 1$.

\medskip

We recall the following definition of weighted sharing:

 \begin {definition}  Let $f$ and $g$ be two non constant meromorphic functions and $k$ be a non-negative integer or $\infty$. For $a \in \mC \cup \{\infty \}$, we denote by $E_k(a, f)$ the set of all a-points of $f$, where an a-point of multiplicity $m$ is counted $m$ times if $m \leq k$ and $k +1$ times if $m > k$. If $E_k(a, f) = E_k(a, g)$, we say that $f$ and $g$ share the value $a$ with weight $k$.
 \end{definition}
 
 We write ``$f$ and $g$ share $(a, k)$" to mean that ``$f$ and $g$ share the value $a$ with weight $k$". Clearly if $f$ and $g$ share $(a, k)$, then $f$ and $g$ share $(a, p)$, $0 \leq p <k$. Also we note that $f$ and $g$ share the value $a$ IM(ignoring multilicity) or CM(counting multiplicity) if and only if $f$ and $g$ share $(a, 0)$ or $(a, \infty)$, respectively.

\begin{definition} let $f$ and $g$ share 1 IM, and let $z_0$ be a zero of $f - 1$ with multiplicity $p$ and a zero of $g - 1$ with multiplicity $q$. We denote by $N^{1)}_E\left(r, 1/(f - 1)\right)$, the counting function of the zeros of $f - 1$ when $p = q = 1$.  By $\overline N^{(2}_E\left(r, 1/(f - 1)\right)$, we denote the counting function of the zeros of $f - 1$ when $p = q \geq 2 $ and by $\overline N_L\left(r, 1/(f - 1)\right)$, we denote the counting function of the zeros of $f - 1$ when $p > q \geq 1$, each point in these counting functions is counted only once; similarly, the terms $N^{1)}_E\left(r, 1/(g - 1)\right)$, $\overline N^{(2}_E\left(r, 1/(g - 1)\right)$ and $\overline N_L\left(r, 1/(g - 1)\right)$. Also, we denote by $\overline N_{f > k}\left(r, 1/(g - 1)\right)$, the reduced counting function of those zeros of $f - 1$ and $g - 1$ such that $p > q = k$, and similarly the term  $\overline N_{g > k}\left(r, 1/(f - 1)\right)$.
\end{definition}

\medskip

With the help of weighted sharing, we generalize Theorem A and Theorem B as 

\begin{theorem}
Let $f$ and $g$ be two non-constant entire functions such that $f$ and $g$ share $0$ CM. Suppose $L(f)$  and $L(g)$ share $(1, l),\;\;l \geq 0$ with one of the following conditions: \\ 
$(i)$ $l \geq 2$ and $\delta(0, f) > 1/2$ \\
$(ii)$ $l = 1$ and $\delta(0, f) > 3/5$\\
$(iii)$ $l = 0$ and $\delta(0, f) > 4/5$.\\
If $\lambda (f) \neq 1$, then $f \equiv g$ unless $L(f).L(g) \equiv 1$.
\end{theorem}

For meromorphic functions, we prove the following result:

\begin{theorem}
Let $f$ and $g$ be two non-constant meromorphic functions of finite order such that $f$ and $g$ share $0$ and $\infty$ CM. Suppose $L(f)$  and $L(g)$ share $(1, l),\;\;l \geq 0$ with one of the following conditions: \\ 
$(i)$ $l \geq 2$ and
\begin{equation}\label{1'}
 (k + 2)\Theta(\infty , f) + 2\delta(0, f) > k + 3
 \end{equation}
$(ii)$ $l = 1$ and 
\begin{equation}\label{2'}
(3k + 5)\Theta(\infty , f) + 5\delta(0, f) > 3k + 9
\end{equation}
$(iii)$ $l = 0$ and 
\begin{equation}\label{3'}
(4k + 5)\Theta(\infty , f) + 5\delta(0, f) > 4k + 9
\end{equation}
If $\lambda (f) \neq 1$ and $\sigma(1/f) \leq \sigma(f)$, then $f \equiv g$ unless $L(f).L(g) \equiv 1$.
\end{theorem}

\medskip

The main tool of our investigations in this paper is Nevanlinna value distribution theory of meromorphic functions(see \cite{HAY}).

\section{Proof of the Main Result}
 We shall use the following results in the proof of our main result: 
 \begin{lemma} \label{a}\cite{LL}
 Let $f$ be a non-constant meromorphic function and $k$ be a non-negative integer. Then
 \begin{equation} \label{1}
 T(r, L(f)) \leq T(r, f) + k\overline N(r, f) + S(r,f).
 \end{equation}
 \end{lemma}
 
 \begin{lemma} \label{b}\cite{LL}
 Let $f$ be a non-constant meromorphic function and $a$ be a meromorphic function such that $T(r, a) = \circ(T(r, f))\;\;\text {as}\;\; r \rightarrow \infty.$ If $f$ is not a polynomial, then
 \begin{align} \label{2}
 N\left(r, \frac{1}{L(f) - L(a)}\right) \leq T(r,L(f)) - T(r, f) + N\left(r, \frac{1}{f - a}\right) + S(r, f)
 \end{align}
 and
 \begin{align}\label{3}
 N\left(r, \frac{1}{L(f)- L(a)}\right) \leq N\left(r, \frac{1}{f - a}\right) + k \overline N(r, f) + S(r, f).
 \end{align}
 \end{lemma}

\begin{lemma} \label{c} \cite{AB}
Let $f$ and $g$ be two non-constant meromorphic functions.\\
{\it (i)} If $f$ and $g$ share $(1, 0)$, then
\begin{align}\label{4}
\overline N_L\left(r, \frac{1}{f-1}\right) \leq \overline N\left(r, \frac{1}{f}\right) + \overline N(r,f) + S(r),
\end{align}
where $S(r) = o(T(r))$ as $r \rightarrow \infty$ with $T(r) = {\text max}\{T(r,f);T(r,g)\}$.\\
{\it (ii)} If $f$ and $g$ share $(1, 1)$, then 
\begin{align}\label{5}
2\overline N_L\left(r, \frac{1}{f-1}\right) + 2 \overline N_L\left(r, \frac{1}{g-1}\right)
& + \overline N^{(2}_E\left(r, \frac{1}{f-1}\right) - \overline N_{f > 2}\left(r, \frac{1}{g-1}\right)\notag\\
& \leq N\left(r, \frac{1}{g - 1}\right) - \overline N\left(r, \frac{1}{g - 1}\right).
\end{align}
\end{lemma}
\begin{lemma} \label{d} \cite{CM}
Suppose $f_j\;(j = 1,2,...,n+1)$ and $g_j\;(j = 1,2,...,n)$ $(n \geq 1)$are entire functions satisfying the following conditions:\\

$(i)$ $\sum_{j = 1}^{n}f_j(z)e^{g_j(z)} \equiv f_{n+1}(z)$,\\

$(ii)$ The order of $f_j(z)$ is less than the order of $e^{g_k(z)}$ for $1 \leq j \leq n+1$, $1 \leq k \leq n$. And furthermore, the order of $f_j(z)$ is less than the order of $e^{g_h(z) - g_k(z)}$ for $n \geq 2$ and $1 \leq j \leq n+1$, $1 \leq h < k \leq m$.\\

Then $f_j \equiv 0 (j = 1,2,...,n+1)$.
\end{lemma}
\begin{lemma}\label{e}\cite{CM}
Suppose $f_j\;(j = 1,2,...,n)$ are meromorphic functions and $g_j\;(j = 1,2,...,n)$ $(n \geq 2)$ are entire functions satisfying the following conditions:\\

$(i)$ $\sum_{j = 1}^{n}f_j(z)e^{g_j(z)} \equiv 0$.\\

$(ii)$ $g_j(z) - g_k(z)$ are non-constants for $1 \leq j < k \leq n$.\\

$(iii)$ For $1 \leq j \leq n$, $1 \leq h < k \leq n$,
$$T(r, f_j) = o(T(r, e^{g_h - g_k})),$$
as $r \rightarrow \infty$. Then $f_j(z) \equiv 0\;\;(j = 1,2,...,n)$.
\end{lemma}

\begin{lemma}\label{f}\cite{CM}
If $h(z)$ be a polynomial of degree $p$ and $f(z) = e^{h(z)}$, then $\lambda(f) = \mu(f) = p$.
\end{lemma}
\begin{lemma}\label{g}\cite{CM}
Let $f(z)$ and $g(z)$ be two non-constant meromorphic functions in the complex plane. If $\lambda(f) < \mu(g)$, then
$T(r, f) = o(T(r, g))$ as $r \rightarrow \infty$.
\end{lemma}
 
\bigskip 

We only prove Theorem 1.4 as the proof of Theorem 1.3 follows on the similar lines.

\bigskip

{\bf Proof of Theorem 1.4:}
First we assume that $L(f) \equiv c$, a finite constant. Then $f$ has to be entire and 
$$f \equiv c_1 + \sum_{i = 1}^m p_i(z) e^{\alpha _{i}z},$$
where $c_1$ is finite constant, $m (\leq k)$ is a positive integer, $\alpha_i$  are distinct complex numbers and $p_i(z)$  are polynomials $(i = 1, 2,...,m)$.\\
Since $\lambda(f) \neq 1$, we get $\lambda(f) < 1$ and so $e^{{\alpha _i}z}$ is constant. Thus $f$ is a polynomial and so $\delta(0, f) = 0$, which contradicts (\ref{1'}), (\ref{2'}) and (\ref{3'}).

\medskip

Assume that both $L(f)$ and $L(g)$ are non-constant. Since $f$ and $g$ share $0$ and $\infty$ CM, and $L(f)$ and $L(g)$ share $(1,l)$, it follows from Milloux's inequality and (\ref{3})
\begin {align*}
T(r,f) & \leq \overline N(r, f) + N(r, \frac{1}{f}) + \overline N\left(r, \frac{1}{L(f) - 1}\right) + S(r, f) \notag\\
& = \overline N(r, g) + N(r, \frac{1}{g}) + \overline N\left(r, \frac{1}{L(g) - 1}\right) + S(r, f) \notag\\
& \leq 2T(r, g) + k\overline N(r, g) + N\left(r, \frac{1}{g - 1}\right) + S(r, f) + S(r, g) \notag\\
& \leq (k + 3)T(r,g) + S(r, f) + S(r, g).
\end{align*}
Similarly
\begin{align*}
T(r,g) \leq (k + 3)T(r,f) + S(r, f) + S(r, g).
\end{align*}
Thus $S(r, f) = S(r, g)$ and $\lambda(f) = \lambda(g)$.

\medskip

Let $F = L(f)$ and $G = L(g)$. Then $F$ and $G$ share $(1, l),\;l \geq 0$. Define
\begin{equation} \label{7}
H = \left(\frac{F''}{F'} - \frac{2F'}{F - 1}\right) - \left(\frac{G''}{G'} - \frac{2G'}{G - 1}\right).
\end{equation}

Assume that $H \not\equiv 0$. Then from (\ref{7}), we have
 $$m(r, H) = S(r, F) + S(r, G).$$
 By the Second fundamental theorem of Nevanlinna, we have
\begin{align} \label{8'}
T(r, F) + T(r, G) & \leq  \overline N(r, F) + \overline N\left(r, \frac{1}{F}\right) + \overline N\left(r, \frac{1}{F - 1}\right) +   \overline N(r, G) + \overline N\left(r, \frac{1}{G}\right)\notag\\
&+ \overline N\left(r, \frac{1}{G - 1}\right) -  N_0\left(r, \frac{1}{F'}\right) - N_0\left(r, \frac{1}{G'}\right) + S(r, F) + S(r, G)\notag\\
 & =  2\overline N(r, F) + \overline N\left(r, \frac{1}{F}\right) + \overline N\left(r, \frac{1}{F - 1}\right) + \overline N\left(r, \frac{1}{G}\right)\notag\\
&+ \overline N\left(r, \frac{1}{G - 1}\right) -  N_0\left(r, \frac{1}{F'}\right) - N_0\left(r, \frac{1}{G'}\right) + S(r, F) + S(r, G),
\end{align}
where  $ N_0(r, 1/F')$ denotes the counting function of the zeros of $F'$ which are not the zeros of $F(F - 1)$ and $ N_0(r, 1/G')$ denotes the counting function of the zeros of $G'$ which are not the zeros of $G(G - 1)$.

\medskip

We consider the following cases:

\medskip

 {\bf Case (i).} If $l \geq 1$, then from (\ref{7}), we have
 \begin{align*}
N^{1)}_E\left(r, \frac{1}{F - 1}\right) & \leq  N \left(r,\frac{1}{H} \right) + S(r, F) + S(r, G) \notag\\
& \leq  T(r, H) + S(r,F) + S(r, G) \notag\\
& =  N(r, H) + S(r,F) + S(r, G) \notag\\
& \leq \overline N_{(2}\left(r, \frac{1}{F}\right) + \overline N_{(2}\left(r, \frac{1}{G}\right) + \overline N_L\left(r, \frac{1}{F - 1}\right)\notag\\
&+ \overline N_L\left(r, \frac{1}{G - 1}\right) +  N_0\left(r, \frac{1}{F'}\right) + N_0\left(r, \frac{1}{G'}\right) + S(r, F) + S(r, G)
\end{align*}
and so
\begin{align}\label{9}
\overline N\left(r, \frac{1}{F - 1}\right)+\overline N\left(r, \frac{1}{G - 1}\right)& = N^{1)}_E\left(r, \frac{1}{F - 1}\right)+\overline N^{(2}_E\left(r, \frac{1}{F - 1}\right) + \overline N_L\left(r, \frac{1}{F - 1}\right)\notag\\
&+ \overline N_L\left(r, \frac{1}{G - 1}\right) + \overline N\left(r, \frac{1}{G - 1}\right) +  S(r, F) + S(r, G) \notag\\
&\leq \overline N_{(2}\left(r, \frac{1}{F}\right) + \overline N_{(2}\left(r, \frac{1}{G}\right) + 2\overline N_L\left(r, \frac{1}{F - 1}\right)\notag\\
&+ 2 \overline N_L\left(r, \frac{1}{G - 1}\right) + \overline N^{(2}_E\left(r, \frac{1}{F - 1}\right) +  \overline N\left(r, \frac{1}{G - 1}\right)\notag\\ 
&+ N_0\left(r, \frac{1}{F'}\right) + N_0\left(r, \frac{1}{G'}\right) +  S(r, F) + S(r, G).
\end{align}
{\bf Subcase 1.1:} When $l = 1$. Then we have
\begin{equation}\label{10}
\overline N_L\left(r, \frac{1}{F - 1}\right) \leq \frac{1}{2}N\left(r, \frac{1}{F'}|F \neq 0\right) \leq \frac{1}{2}\overline N(r, F) + \frac{1}{2}\overline N\left(r, \frac{1}{F}\right),
\end{equation}
where $N\left(r, \frac{1}{F'}|F \neq 0\right)$ denotes the zeros of $F'$, that are not the zeros of $F$.

\medskip

 From (\ref{5}) and (\ref{10}), we have
\begin{align} \label{11}
2\overline N_L\left(r, \frac{1}{F - 1}\right) + 2 \overline N_L\left(r, \frac{1}{G - 1}\right)& + \overline N^{(2}_E\left(r, \frac{1}{F - 1}\right) + \overline N\left(r, \frac{1}{G - 1}\right)\notag\\
&\leq N\left(r, \frac{1}{G - 1}\right) + \overline N_L\left(r, \frac{1}{F - 1}\right) + S(r, F) + S(r, G)\notag \\
&\leq N\left(r, \frac{1}{G - 1}\right) + \frac{1}{2}\overline N(r, F) + \frac{1}{2}\overline N\left(r, \frac{1}{F}\right) + S(r, F) + S(r, G).
\end{align}
Thus, from (\ref{9}) and (\ref{11}), we have
\begin{align} \label{12}
\overline N\left(r, \frac{1}{F - 1}\right)+\overline N\left(r, \frac{1}{G - 1}\right) & \leq  \frac{1}{2}\overline N(r, F) + \overline N_{(2}\left(r, \frac{1}{F}\right) + \overline N_{(2}\left(r, \frac{1}{G}\right)\notag\\
&+ \frac{1}{2}\overline N\left(r, \frac{1}{F}\right) + N\left(r, \frac{1}{G - 1}\right)\notag \\
&+ N_0\left(r, \frac{1}{F'}\right) + N_0\left(r, \frac{1}{G'}\right) + S(r, F) + S(r, G)\notag\\
& \leq \frac{1}{2}\overline N(r, F) + \overline N_{(2}\left(r, \frac{1}{F}\right) + \overline N_{(2}\left(r, \frac{1}{G}\right)\notag\\
&+\frac{1}{2}\overline N\left(r, \frac{1}{F}\right) + T(r, G)\notag \\
&+ N_0\left(r, \frac{1}{F'}\right) + N_0\left(r, \frac{1}{G'}\right) + S(r, F) + S(r, G).
\end{align}
From (\ref{2}), (\ref{3}), (\ref{8'}) and (\ref{12}), we obtain
\begin{align*}
T(r,F) & \leq \frac{5}{2}\overline N(r, F) + \overline N\left(r, \frac{1}{F}\right) + \overline N_{(2}\left(r, \frac{1}{F}\right) +  \overline N\left(r, \frac{1}{G}\right)+
\overline N_{(2}\left(r, \frac{1}{G}\right)\notag\\
& + \frac{1}{2}\overline N\left(r, \frac{1}{F}\right)+ S(r, F) + S(r, G)\notag\\
& \leq  \frac{5}{2}\overline N(r, F) + N\left(r, \frac{1}{F}\right) + N\left(r, \frac{1}{G}\right) + \frac{1}{2}\overline N\left(r, \frac{1}{F}\right)+ S(r, F) + S(r, G)\notag\\
& = \frac{5}{2}\overline N(r, F) + N\left(r, \frac{1}{L(f)}\right) + \frac{1}{2}\overline N\left(r, \frac{1}{L(f)}\right) + N\left(r, \frac{1}{L(g)}\right) + S(r, F) + S(r, G)\notag\\
& \leq \frac{5}{2}\overline N(r, f) + T(r, L(f)) - T(r,f) +  N\left(r, \frac{1}{f}\right) + \frac{1}{2}N\left(r, \frac{1}{f}\right)\notag\\
& + \frac{k}{2}\overline N(r, f)+ N(r, \frac{1}{g}) + k\overline N(r, g)+ S(r, f) + S(r, g)\notag\\
& = T(r, L(f)) - T(r,f) + \left(\frac{3k + 5}{2}\right)\overline N(r, f) + \frac{5}{2} N\left(r, \frac{1}{f}\right) + S(r, f).
\end{align*}
That is,
$$2T(r, f) \leq  (3k + 5)\overline N(r, f) + 5 N\left(r, \frac{1}{f}\right) + S(r, f),$$
and so $(3k + 5)\Theta(\infty , f) + 5 \delta(0, f) \leq 3k + 8$, a contradiction to (\ref{2'}).
\medskip

{\bf Subcase 1.2:} When $l \geq 2$.\\
In this case, we have  

\medskip

$2\overline N_L\left(r, \frac{1}{F - 1}\right) + 2 \overline N_L\left(r, \frac{1}{G - 1}\right) + \overline N^{(2}_E\left(r, \frac{1}{F - 1}\right) + \overline N\left(r, \frac{1}{G - 1}\right)$
 $$\leq N\left(r, \frac{1}{G - 1}\right) + S(r, F) + S(r, G).$$
Thus from (\ref{9}), we get
\begin{align} \label{13}
\overline N\left(r, \frac{1}{F - 1}\right)+\overline N\left(r, \frac{1}{G - 1}\right)& \leq  \overline N_{(2}\left(r, \frac{1}{F}\right) + \overline N_{(2}\left(r, \frac{1}{G}\right) +  N\left(r, \frac{1}{G - 1}\right)\notag\\
&+ N_0\left(r, \frac{1}{F'}\right) + N_0\left(r, \frac{1}{G'}\right) + S(r, F) + S(r, G)\notag\\
& \leq \overline N_{(2}\left(r, \frac{1}{F}\right) + \overline N_{(2}\left(r, \frac{1}{G}\right) +  T(r, G)\notag\\
&+ N_0\left(r, \frac{1}{F'}\right) + N_0\left(r, \frac{1}{G'}\right) + S(r, F) + S(r, G).
\end{align}
Since $f$ and $g$ share $0$ and $\infty$ CM, from (\ref{2}), (\ref{3}), (\ref{8'}) and (\ref{13}), we obtain
\begin{align*}
T(r,F) &\leq  2\overline N(r, F) + \overline N\left(r, \frac{1}{F}\right) + \overline N_{(2}\left(r, \frac{1}{F}\right) + \overline N\left(r, \frac{1}{G}\right) + \overline N_{(2}\left(r, \frac{1}{G}\right) + S(r, F) + S(r, G) \notag\\
& \leq 2\overline N(r, f) + N\left(r, \frac{1}{F}\right) + N\left(r, \frac{1}{G}\right) + S(r, F) + S(r, G) \notag\\
& =  2\overline N(r, f) + N\left(r, \frac{1}{L(f)}\right) + N\left(r, \frac{1}{L(g)}\right) + S(r, f) + S(r, g)\notag\\
& \leq 2\overline N(r, f) + T(r, L(f)) - T(r,f) +  N\left(r, \frac{1}{f}\right) + N(r, \frac{1}{g}) + k\overline N(r, g) + S(r, f) + S(r, g)\notag\\
& = T(r, L(f)) - T(r,f) + (k + 2)\overline N(r, f) + 2N\left(r, \frac{1}{f}\right) + S(r, f).
\end{align*}
 That is,
$$T(r, f) \leq  (k + 2)\overline N(r, f) + 2N\left(r, \frac{1}{f}\right) + S(r, f),$$
and so $(k + 2)\Theta(\infty , f) + 2\delta(0, f) \leq k + 3$, a contradiction to (\ref{1'}).

{\bf Case (ii).} If $l = 0$, then we have
 
 \medskip
 
$N^{1)}_E\left(r, \frac{1}{F - 1}\right) = N^{1)}_E\left(r, \frac{1}{G - 1}\right) + S(r, F) + S(r, G)$,

\medskip

$\overline N^{(2}_E\left(r, \frac{1}{F - 1}\right) = \overline N^{(2}_E\left(r, \frac{1}{G - 1}\right) + S(r, F) + S(r, G)$,

\medskip

and also from (\ref{7}), we have
\begin{align}\label{14}
\overline N\left(r, \frac{1}{F - 1}\right)+\overline N\left(r, \frac{1}{G - 1}\right)&\leq N^{1)}_E\left(r, \frac{1}{F - 1}\right)+\overline N^{(2}_E\left(r, \frac{1}{F - 1}\right) + \overline N_L\left(r, \frac{1}{F - 1}\right)\notag\\
&+ \overline N_L\left(r, \frac{1}{G - 1}\right) + \overline N\left(r, \frac{1}{G - 1}\right) + S(r, F) + S(r, G) \notag\\
& \leq  N^{1)}_E\left(r, \frac{1}{F - 1}\right) + \overline N_L\left(r, \frac{1}{F - 1}\right) + N\left(r, \frac{1}{G - 1}\right)\notag\\
& + S(r, F) + S(r, G) \notag\\
& \leq \overline N_{(2}\left(r, \frac{1}{F}\right) + \overline N_{(2}\left(r, \frac{1}{G}\right) + 2\overline N_L\left(r, \frac{1}{F - 1}\right)\notag\\
&+ \overline N_L\left(r, \frac{1}{G - 1}\right) + N\left(r, \frac{1}{G - 1}\right) + N_0\left(r, \frac{1}{F'}\right)\notag\\
&+ N_0\left(r, \frac{1}{G'}\right) + S(r, F) + S(r, G).
\end{align}
 From (\ref{2}), (\ref{3}), (\ref{4}), (\ref{8'}) and (\ref{14}), we obtain
\begin{align*}
T(r,F) &\leq 2 \overline N(r, F) + \overline N\left(r, \frac{1}{F}\right) + \overline N_{(2}\left(r, \frac{1}{F}\right) + \overline N\left(r, \frac{1}{G}\right) + \overline N_{(2}\left(r, \frac{1}{G}\right)\notag\\
&+ 2\overline N_L\left(r, \frac{1}{F - 1}\right) + \overline N_L\left(r, \frac{1}{G - 1}\right) + S(r, F) + S(r, G) \notag\\
& \leq  2 \overline N(r, F) + N\left(r, \frac{1}{F}\right) + N\left(r, \frac{1}{G}\right) + 2\overline N\left(r, \frac{1}{F}\right) + 2\overline N(r, F)\notag\\
&+ \overline N\left(r, \frac{1}{G}\right) + \overline N(r, G) + S(r, F) + S(r, G) \notag\\
& \leq 5 \overline N(r, f) + N\left(r, \frac{1}{L(f)}\right) + 2N\left(r, \frac{1}{L(f)}\right) + 2N\left(r, \frac{1}{L(g)}\right) + S(r, F) + S(r, G) \notag\\
& \leq 5 \overline N(r, f) + T(r, L(f)) - T(r,f) +  N\left(r, \frac{1}{f}\right) +  2N\left(r, \frac{1}{f}\right) + 2k \overline N(r, f)\notag\\
& + 2N(r, \frac{1}{g}) + 2k\overline N(r, g) + S(r, f) + S(r, g)\notag\\
& \leq T(r, L(f)) - T(r,f) + (4k + 5)\overline N(r, f) + 5 N\left(r, \frac{1}{f}\right) + S(r, f).
\end{align*}
 That is,
$$T(r, f) \leq  (4k + 5)\overline N(r, f) + 5N\left(r, \frac{1}{f}\right) + S(r, f),$$
and so $(4k + 5)\Theta(\infty , f) + 5\delta(0, f) \leq 4k + 9$, a contradiction to (\ref{3'}).
\medskip

Thus our supposition is wrong and hence $H \equiv 0$. So (\ref{7}) implies that
$$\frac{F''}{F'} - \frac{2F'}{F - 1} = \frac{G''}{G'} - \frac{2G'}{G - 1},$$
and so we obtain
\begin{equation}\label{15}
\frac{1}{F - 1} = \frac{C}{G - 1} + D,
\end{equation} 
where $C \neq 0$ and $D$ are constants.\\

Here, the following three cases can arise:

\medskip

{\bf Case$(a):$} When $D \neq 0,\;-1$. We rewrite (\ref{15}) as
$$\frac{G - 1}{C} = \frac{F - 1}{D + 1 - DF},$$
we have
$$\overline N(r, G) = \overline N\left(r,\frac{1}{F - (D + 1)/D}\right).$$

\medskip

By Second fundamental theorem of Nevanlinna and (\ref{2}), we have
\begin{align*}
T(r, L(f)) & =  T(r, F) + S(r,f)\\
& \leq  \overline N(r,F) + \overline N\left(r,\frac{1}{F}\right) + \overline N\left(r,\frac{1}{F - (D + 1)/D}\right) + S(r,f)\\
& \leq \overline N(r,F) + \overline N\left(r,\frac{1}{F}\right) + \overline N(r, G) + S(r,f)\\
& \leq N\left(r,\frac{1}{L(f)}\right) + 2\overline N(r, f) + S(r,f)\\
& \leq T(r, L(f)) - T(r, f) + 2\overline N(r, f) + N(r, \frac{1}{f}) + S(r, f).
\end{align*}
Thus 
$$T(r, f) \leq 2\overline N(r, f) + N(r, \frac{1}{f}) + S(r, f),$$
and so $2 \Theta(\infty , f) + \delta(0, f) \leq 2,$ which contradicts (\ref{1'}),(\ref{2'}) and (\ref{3'}).

\medskip 

{\bf Case$(b):$} When $D = 0$. Then from (\ref{15}), we have
\begin{equation}\label{16}
G = CF - (C - 1).
\end{equation}
So if $C \neq 1$, then
$$\overline N\left(r, \frac{1}{G}\right) = \overline N\left(r,\frac{1}{F - (C - 1)/C}\right).$$

\medskip

Since $f$ and $g$ share $0$ and $\infty$ CM, by Second fundamental theorem of Nevanlinna, (\ref{2}) and (\ref{3}) gives
\begin{align*}
T(r, L(f)) & =  T(r, F) + S(r,f)\\
& \leq  \overline N(r,F) + \overline N\left(r,\frac{1}{F}\right) + \overline N\left(r,\frac{1}{F - (C - 1)/C}\right) + S(r,f)\\
& =  \overline N(r, f) + \overline N\left(r,\frac{1}{F}\right) + \overline N(r, \frac{1}{G}) + S(r,f)\\
& \leq \overline N(r, f) +  N\left(r,\frac{1}{L(f)}\right) + N\left(r,\frac{1}{L(g)}\right) + S(r,f)\\
& \leq \overline N(r, f) + T(r, L(f)) - T(r, f) + N(r, \frac{1}{f}) + N(r, \frac{1}{g}) + k\overline N(r, g) +  + S(r, f)\\
& = T(r, L(f)) - T(r, f) + (k + 1)\overline N(r, f) + 2 N(r, \frac{1}{f}) + S(r, f).
\end{align*}
Thus 
$$T(r, f) \leq (k + 1)\overline N(r, f) +2 N(r, \frac{1}{f}) + S(r, f),$$
and so $(k + 1) \Theta(\infty , f) + 2 \delta(0, f) \leq k + 2,$ which contradicts (\ref{1'}),(\ref{2'}) and (\ref{3'}).

\medskip

Thus, $C = 1$ and so in this case from (\ref{16}), we obtain $F \equiv G$ and so
$$L(f) \equiv L(g).$$

{\bf Case$(c):$} When $D = -1$. Then from (\ref{15}) we have
\begin{equation} \label{17}
\frac{1}{F - 1} = \frac{C}{G - 1} - 1.
\end{equation}
So if $C \neq -1$, then
$$\overline N\left(r, \frac{1}{G}\right) = \overline N\left(r,\frac{1}{F - C/(C + 1)}\right).$$
Since $f$ and $g$ share $0$ and $\infty$ CM, by Second fundamental theorem of Nevanlinna, (\ref{2}) and (\ref{3}), we have
\begin{align*}
T(r, L(f)) & =  T(r, F) + S(r,f)\\
& \leq  \overline N(r,F) + \overline N\left(r,\frac{1}{F}\right) + \overline N\left(r,\frac{1}{F - C/(C+1)}\right) + S(r,f)\\
& =  \overline N(r, f) + \overline N\left(r,\frac{1}{F}\right) + \overline N(r, \frac{1}{G}) + S(r,f)\\
& \leq \overline N(r, f) +  N\left(r,\frac{1}{L(f)}\right) + N\left(r,\frac{1}{L(g)}\right) + S(r,f)\\
& \leq \overline N(r, f) + T(r, L(f)) - T(r, f) + N(r, \frac{1}{f}) + N(r, \frac{1}{g}) + k\overline N(r, g) + S(r, f)\\
& = T(r, L(f)) - T(r, f) + (k + 1)\overline N(r, f) + 2 N(r, \frac{1}{f}) + S(r, f).
\end{align*}
Thus 
$$T(r, f) \leq (k + 1)\overline N(r, f) +2 N(r, \frac{1}{f}) + S(r, f),$$
and so $(k + 1) \Theta(\infty , f) + 2 \delta(0, f) \leq k + 2,$ which contradicts (\ref{1'}),(\ref{2'}) and (\ref{3'}).
\medskip

Thus, $C = -1$ and so in this case from (\ref{17}), we obtain $FG \equiv 1$ and so $L(f)L[f] = 1.$ 

\medskip

If $L(f) \equiv L(g)$, then $L(f - g) \equiv 0$ and so $f - g$ has to be entire and we have (see \cite{LA}) 
$$f - g = \sum_{i = 1}^m p_i(z) e^{\alpha _{i}z},$$
where $m (\leq k)$ is a positive integer, $\alpha_i$  are distinct complex numbers and $p_i(z)$  are polynomials $(i = 1, 2,...,m)$.\\
Thus
 $$\lambda(f - g) = \lambda \left(\sum_{i = 1}^m p_i(z) e^{\alpha _{i }z}\right) \leq 1.$$
 We consider the following cases:
 
 \medskip
 
 {\bf Case (i).} When $\lambda(f) < 1$. Since $f$ and $g$ share $0$ and $\infty$ CM, we have $f/g = e^{h(z)}$, where $h(z)$ is an entire function. Also as $\lambda(f) = \lambda(g)$, we have
$$\lambda(e^{h(z)}) = \lambda(f/g) \leq {\text max}\{\lambda(f), \lambda(1/g)\} < 1.$$
Thus $e^{h(z)}$ is a constant, say $c$ and so $f \equiv c g$ which implies that $L(f) \equiv c L(g)$. But $L(f) \equiv L(g)$, so we get $c = 1$ and thus $f \equiv g$.

\medskip

{\bf Case (ii).} When $\lambda(f) > 1$. Since $f$ and $g$ are meromorphic functions of finite order, by Hadamard's factorization theorem we have
$$f(z) = \frac{P(z)}{Q(z)}e^{l_1(z)} \;\; {\text and}\;\;g(z) = \frac{P(z)}{Q(z)}e^{l_2(z)},$$
where $P(z)$ is the canonical product formed with the common zeros of $f$ and $g$, $Q(z)$ is the canonical product formed with the common poles of $f$ and $g$, and $l_1,\;l_2$ are the polynomials of degree less than or equal to $\lambda(f),\;\lambda(g)$ repectively.\\
Thus 
$$f - g = \frac{P(z)}{Q(z)}e^{l_1(z)} - \frac{P(z)}{Q(z)}e^{l_2(z)},$$
 or we can write
\begin{equation}
\frac{P(z)}{Q(z)}e^{l_1(z)} - \frac{P(z)}{Q(z)}e^{l_2(z)} - (f - g)e^{l_3(z)} \equiv 0,
\end{equation} 
where $l_3(z) \equiv 0$.

\medskip

Also 
$$\lambda(P) = \sigma(f) \leq \sigma(f - g) \leq \lambda(f - g) \leq 1,$$
and since $\sigma(1/f) \leq \sigma(f)$, we have
$$\lambda(Q) = \sigma(1/f) \leq \sigma(f) \leq \sigma(f - g) \leq \lambda(f - g) \leq 1.$$
Thus 
$$\lambda\left(\frac{P}{Q}\right) \leq\; {\text max}\{\lambda(P), \lambda(Q)\} \leq 1.$$
Since $f - g = (e^{l_1 - l_2})g$ and $\lambda(f) = \lambda(g) > 1$, we have $\lambda(e^{l_1}) > 1$, $\lambda(e^{l_1}) > 1$ and $\lambda(e^{l_1 - l_2}) > 1$, and so $\lambda(e^{l_i - l_j}) > 1$, where $1\leq i < j \leq 3$. Thus $l_i - l_j$ is non-constant, where  $1\leq i < j \leq 3$ and by lemma \ref{f} and \ref{g}, we get
$$T(r, f - g) = o(T(r, e^{l_i - l_j}))\;\;{\text and}\;\;T(r, P/Q) = o(T(r, e^{l_i - l_j})),$$
as $r \rightarrow \infty$.
Thus by lemma(\ref{e}), we have  $P/Q \equiv 0$ and $f - g \equiv 0$ which implies that $f(z) \equiv 0$, which is a contradiction. So $l_1 = l_2$ and hence $f \equiv g$.

~~~~~~~~~~~~~~~~~~~~~~~~~~~~~~~~~~~~~~~~~~~~~~~~~~~~~~~~~~~~~~~~~~~~~~~~~~~~~~~~~~~~~~~~~~~~~~~~~~~~~$\Box$

\bibliographystyle{amsplain}

\end{document}